\documentclass[leqno,11pt]{amsart}

\newtheorem{Thm}{Theorem}

\newtheorem{Rem}{Remark}
\newtheorem{Prop}{Proposition}
\newtheorem{Cor}{Corollary}
\newtheorem{Lem}{Lemma}

\def\longformule#1#2{
\displaylines{
\qquad{#1}
\hfill\cr
\hfill {#2}
\qquad\cr
}
}

\def \R{{\mathbf R}}

\def \eps{{\varepsilon}}
\def \e{{\varepsilon}}

\def \RR{{\mathcal R}}

\def \d{{\partial}}

\def \GRAD{\nabla\!}
\def \DIV{\nabla\! \! \cdot }

\def\eqdefa{\buildrel\hbox{\footnotesize def}\over =}

\newcommand{\dps}{\displaystyle}

\begin{document}
\title[Pressureless magneto--hydrodynamics]{
Asymptotic results for pressureless magneto--hydrodynamics}

\author[I. Gallagher]{Isabelle Gallagher}
\address[I. Gallagher]%
{Centre de Math{\'e}matiques UMR 7640 \\
     Ecole Polytechnique \\
91128 Palaiseau
    \\
     FRANCE}
\email{Isabelle.Gallagher@math.polytechnique.fr}

\author[L. Saint-Raymond]{Laure Saint-Raymond}
\address[L. Saint-Raymond]%
{ Laboratoire J.-L. Lions UMR 7598\\ Universit{\'e} Paris
VI\\
175, rue du Chevaleret\\ 75013 Paris\\FRANCE }
\email{saintray@ann.jussieu.fr }
\keywords{rotating pressureless gas, asymptotic behaviour, oscillations}
\subjclass[2000]{Primary 35B40; Secondary 76U05,76W05}

\maketitle

{\small{\bf Abstract.}
We are interested in the life span and the asymptotic behaviour of the
solutions to a system governing the motion of a pressureless gas,
submitted to a strong, inhomogeneous magnetic field~$ \e^{-1} B(x)$,
of variable amplitude but fixed direction --- this is a first step in
the direction of the study of rotating Euler equations.  This leads to the
study of
a multi--dimensional Burgers type system on the velocity field~$
u_\e$, penalized by a rotating term~$  \e^{-1} u_\e \wedge B(x)$. We
prove that the unique, smooth solution of this Burgers  system
exists on a uniform time interval~$ [0,T]$. We also prove that the phase
of oscillation of~$ u_\e$ is an order one perturbation of the phase
obtained in the case of a pure rotation (with no nonlinear transport
term), $ \e^{-1}B(x)t$. Finally going back to the pressureless gas
system, we obtain the asymptotics of the density as~$ \e$ goes to zero.}

\vskip .4cm\noindent 

\centerline{\bf R\'esultats asymptotiques pour la
  magn\'eto-hydrodynamique sans pression}

\vskip .4cm\noindent 

{\small{\bf R\'esum\'e.}
 On s'int\'eresse au temps d'existence et au
  comportement asymptotique des solutions d'un syst\`eme mod\'elisant
  le mouvement d'un gaz sans pression, soumis \`a un fort champ
  magn\'etique~$ \e^{-1} B(x)$, d'intensit\'e variable mais de
  direction fixe --- cela \'etant un premier pas dans la
  compr\'ehension des \'equations d'Euler en rotation rapide. Cela
  conduit \`a l'\'etude d'un syst\`eme de type Burgers
  multi--dimensionnel sur le champ de vitesse~$
u_\e$, p\'enalis\'e  par un terme de rotation~$  \e^{-1} u_\e \wedge
  B(x)$. On d\'emontre que la solution unique r\'eguli\`ere de ce
  syst\`eme de Burgers existe sur un intervalle de temps uniforme~$
  [0,T]$. On montre aussi que la phase d'oscillation de~$ u_\e$ est
  une perturbation au premier ordre de la phase obtenue dans le cas
  d'une rotation pure (sans terme de transport non lin\'eaire), $
  \e^{-1}B(x)t$. Enfin en revenant au syst\`eme des gaz sans pression,
  on obtient le comportement asymptotique de la densit\'e quand~$ \e$
  tend vers z\'ero.}

\section{Introduction}

The aim of this paper is to study the asymptotic behaviour of a fluid
submitted to a strong external inhomogeneous magnetic field.

The case when the field is
constant has been studied by a number of
authors, both for compressible and incompressible models of fluids
(see for instance \cite {BMN}, \cite{cdgg} or~\cite{Grenier}
 for incompressible fluids, and
\cite{FS} or
\cite{GSR} for rarefied plasmas). In that case, one can not only derive the
asymptotic average motion (which is given by the weak limit  of the velocity
field), but one can also describe all the oscillations in the system and
possibly their coupling~: the filtering techniques used for that rely on
explicit computations in  Fourier space.

In the case when the magnetic field is inhomogeneous, those methods are not
relevant any more.  Weak compactness and compensated compactness arguments
allow nevertheless to determine the average motion (see \cite{GSR} in the case
of a rarefied plasma governed by the Vlasov-Poisson system, and \cite{GaSR} in
the case of a viscous incompressible fluid). In order to describe the
oscillating component of the motion, one has to understand the interaction
between the penalization and the nonlinear term of transport~: indeed one
expects that the flow modifies substantially the phase of oscillation (which
is of course inhomogeneous).

We propose here to analyse this interaction for a simplified model of
magneto-hydro-dynamics, the so-called  Euler system of  pressureless gas
dynamics.

\subsection{A simple model for magneto-hydro-dynamics}

We  consider the following system of partial differential equations~:
\begin{equation}
\label{Euler3D}
\begin{aligned}
\d _t \rho +\DIV (\rho u)=0,\quad x\in \R^3, t>0\\
 \d _t (\rho u)+\DIV (\rho u\otimes u)= \rho u\wedge B,\quad x\in \R^3, t>0\\
\rho(t=0)\equiv \rho_0,\quad u(t=0)\equiv u_0\quad x\in \R^3,
\end{aligned}
\end{equation}
where $\rho$ denotes the density of the fluid, $u$ its mean velocity and $B$
is the external magnetic field ($\DIV B=0$). The first equation expresses the
local conservation of mass, while the second one gives the local conservation
of momentum provided that there is no internal force (no pressure). This
assumption is relevant only in some particular regimes (corresponding to sticky
particles \cite{BrGr}). From a physical point of view, this may seem a strong
restriction, but it allows to perform a first mathematical study of that
type of inhomogeneous singular perturbation problem~: indeed in this special
case a major simplification arises since the equation on the mean velocity
can be (at least formally) decoupled from the rest of the system~:
$$
\d _t u+(u\cdot \GRAD )u =u\wedge B\quad x\in \R^3, t>0.
$$
We then obtain a system of Burgers' type, that is a prototype of hyperbolic
system. A work in progress should extend the present results to more realistic
models, in particular to the 3D incompressible Euler system.

\bigskip

In order to further simplify the analysis, we assume that the direction of
the field
$B$ is constant
$$B(x)\equiv \frac1\eps b(x_1,x_2) e_3,\quad (x_1,x_2)\in \R^2\, \hbox{ and
} e_3= {}^t(0,0,1)$$  which allows to get rid of
the geometry of the field lines (for detailed comments on this subject see  for
instance
\cite{GaSR}, Remark 1.4). Any solution to the system (\ref{Euler3D})
has then  uniform regularity with respect to the variable $x_3$. To
isolate the phenomenon of inhomogeneous oscillations with instantaneous loss
of regularity, we restrict therefore our attention to the 2D singular
perturbation problem, in the plane orthogonal to the magnetic field. We
finally have~:
\begin{equation}
\label{Euler}
\begin{aligned}
\d _t \rho +\DIV (\rho u)=0,\quad x\in \R^2, t>0\\
 \d _t (\rho u)+\DIV (\rho u\otimes u)=\frac b\eps \rho u^\perp,\quad x\in
\R^2, t>0\\
\rho(t=0)\equiv \rho_0,\quad u(t=0)\equiv u_0,\quad x\in \R^2,
\end{aligned}
\end{equation}
where $u^\perp$ denotes the vector field with components $(u_2, -u_1)$, and the
intensity $b$ of the magnetic field satisfies the following assumptions~:
$$b \in C^\infty(\R^2)\cap W^{2,\infty}(\R^2), \eqno(H0)$$
$$ \inf _{x\in \R^2} b(x) =b_->0. \eqno (H1)$$

\bigskip

A standard fixed point argument then allows  to prove the local well-posedness
of (\ref{Euler}). The result is the following.

\begin{Thm}
\label{thm:existenceeps}
{\sl
Consider a  function $b$ satisfying assumptions $(H0)(H1)$.
Let $\rho_0$  be a nonnegative function  and $u_0$ be a
vector-field in $H^s(\R^2)$ ($s>2$).
Then, for all $\eps>0$, there
exist   $T_\eps\in ]0,+\infty]$ and a unique solution  of~(\ref{Euler}),  $(\rho_\eps,u_\eps) \in
L^\infty_{loc}([0,T_\eps[,
H^s(\R^2))$.
}
\end{Thm}

\noindent
Note that the lifespan $T_\eps$ of the solution depends on $\eps$, and that the
lower bound on $T_\eps$ coming from the Duhamel formula goes to zero as
$\eps \to 0$. The first difficulty to study the asymptotics $\eps \to 0$
consists then in understanding why the magnetic penalization does not
destabilize the system, and in proving that the solution $(\rho_\eps,u_\eps)$
exists on a uniform interval of time.

\subsection{Formal analysis}

Before stating more precise results on the lifespan of the solutions 	and on
the asymptotics $\eps \to 0$, we have chosen to give some simple observations
about the problem to guide  intuition. In this first approach we restrict
our attention to the analysis of the equation governing the  velocity.

\bigskip

The first step of the formal analysis consists in determining the mean
behaviour of the velocity field, that is the weak limit of $u_\eps$. We have
$$u_\eps=\frac\eps b \left( \d _t u_\eps +(u_\eps\cdot \GRAD
)u_\eps\right)^\perp\,.$$
As $b$ is bounded from below, if we are able to establish convenient a priori
bounds on $u_\eps$, this will imply
$$u_\eps \rightharpoonup 0$$
in some weak sense. This means that we expect the velocity to oscillate at
high frequency (on vanishing temporal or spatial  scales).

\bigskip
Another  way to get an idea of the asymptotic behaviour of the velocity
is to study the simple case when $b$ is constant. The group of oscillations
generated by the magnetic penalization is then homogeneous~:
$$R\left( \frac t\eps\right) u= u\cos \left( \frac {bt}\eps\right) -u^\perp
\sin  \left( \frac {bt}\eps\right)\,,$$
which corresponds to the rotation with frequency $2\pi b/\eps$. As the
coefficients are constant, this group is not perturbed by the transport.
Classical filtering methods (see namely \cite{Grenier},\cite{schochet}) can
then be applied~: setting
$$v_\eps\eqdefa R\left( -\frac t\eps\right) u_\eps$$
leads to
$$\d _t v_\eps+ Q\left(\frac t\eps, v_\eps,v_\eps\right)=0$$
where $Q\left(\frac t\eps, .,.\right)$ is a quadratic form with bounded
coefficients depending on $t/\eps$. As there is only one oscillation
frequency, 
there is no resonance, which implies that
$$v_\eps \to u_0$$
in some strong sense, provided that convenient a priori bounds on $u_\eps$ (and
consequently on
$v_\eps$) hold. This means that we can describe completely the oscillations
and get a strong convergence result. Of course, we get as a corollary that the
lifespan $T_\eps$ is uniformly bounded from below, and we even expect that
$T_\eps \to +\infty$ as $\eps \to 0$.

\bigskip

The case we consider here is much more complicated. The group of oscillations
generated by the magnetic penalization is again very easy to describe~:
$$R\left( \frac t\eps,x\right) u= u\cos \left( \frac {b(x)t}\eps\right)
-u^\perp
\sin  \left( \frac {b(x)t}\eps\right)\,,$$
but it is non homogeneous, which entails

\noindent
$\bullet$ a loss of regularity ($R\left( \frac t\eps,x\right) u$ blows up in
all Sobolev norms $H^s(\R^2)$ for $s>0$);

\noindent
$\bullet$ an interaction with the transport operator (with the same definition
of $v_\eps=R\left( -\frac t\eps,x\right) u_\eps$ as previously, we do not
expect $\d_t v_\eps$ to be bounded in any space of distributions).

\noindent
The stake behind this model problem is to understand how to overcome these
difficulties. The first step is to explain how the phase of oscillations is
modified by the flow~: note that even a small correction on the phase changes
strongly the vector field. Then we have to establish  a strong convergence
result using a new method~: classical energy methods fail because of the lack
of regularity on approximate solutions.
Here an appropriate rewriting of the system  by means of
characteristics associated with the flow allows to understand the underlying
structure and to answer both questions~: in particular we will see that the
spaces which are well adapted for this type of study are constructed on
$L^\infty(\R^2)$. In the case of  incompressible dynamics, the analysis will
be therefore much more difficult since the transport is replaced by a
non-local pseudo-differential operator.

\subsection{Main results}

As long as the solution $(\rho_\eps,u_\eps)$ of system (\ref{Euler}) is
smooth, the velocity $u_\eps$ satisfies the equation of Burgers type
\begin{equation}
\label{Burgers}
\begin{aligned}\d _t u_\eps+(u_\eps\cdot \GRAD) u_\eps+\frac b\eps
u_\eps^\perp=0,\quad x\in \R^2,t>0\\
u_\eps(t=0)=u^0,\quad x\in \R^2.
\end{aligned}
\end{equation}
Using refined a priori estimates on this last equation, we can prove that for
all $\eps >0$ it admits a smooth solution  on a uniform time $T>0$. We
will prove the following result. 

\begin{Thm}\label{existence}
{\sl
Consider a  function $b$ satisfying assumptions $(H0)(H1)$.
Let $\rho_0$  be a nonnegative function in
$W^{s-1,\infty}(\R^2)$, and $u_0$ be  a vector-field in
$W^{s,\infty}(\R^2)$ ($s\geq 1$).
Then,  there
exists   $T^*\in ]0,+\infty]$ such that, for all $T<T^*$ and all $\eps \leq
\eps_T$, there is  a unique
$(\rho_\eps,u_\eps)
\in  L^\infty([0,T],
W^{s-1,\infty}(\R^2)\times W^{s,\infty}(\R^2))$ solution  of (\ref{Euler})
(which is nevertheless not uniformly bounded in
$L^\infty([0,T],
W^{s-1,\infty}(\R^2)\times W^{s,\infty}(\R^2))$ for
$s>0$).
}
\end{Thm}

In this framework, it is relevant to consider the asymptotics $\eps \to 0$ on
the time interval $[0,T]$. The same type of computations as previously allows
to prove that the velocity field behaves almost as in the constant case (with
slight modifications of the phase of oscillations).

\begin{Thm}\label{convergence-u}
{\sl
Consider a  function $b$ satisfying assumptions $(H0)(H1)$.
Let  $u_0$ be  a vector-field in
$W^{s,\infty}(\R^2)$ ($s\geq 1$). For all $T\leq T^*$ as in Theorem
\ref{existence} and all $\eps\leq
\eps_T$, denote by
$u_\eps$ the solution of (\ref{Burgers}) in $L^\infty([0,T],
W^{s,\infty}(\R^2))$. Then,
$$u_\eps(t,x)- \left( u_0(x) \cos \theta_\eps(t,x)-u_0^\perp(x) \sin
\theta_\eps(t,x)
\right)$$
converges strongly to 0 in $L^\infty([0,T]\times \R^2)$, where the phase
$\theta_\eps$ is defined by the following equation
\begin{equation}
\label{theta}
\theta_\e (t,x) = \frac{b(x)t}{\e} - t u_0(x) \cdot \nabla \log b(x)
\sin \theta_\e (t,x)
\end{equation}
$$ \quad \quad \quad \quad \quad \quad \quad \quad \quad \quad \quad \quad 
+ t
u_0^\perp(x)\cdot \nabla \log b(x) \cos \theta_\e (t,x).
$$
}
\end{Thm}

Rewriting the equation on the density $\rho_\eps$ with a transport term
and a penalization term (coming from the divergence of $u_\eps$ which is of
order $1/\eps$)
\begin{equation}
\label{density}
\begin{aligned}\d _t \rho_\eps+(u_\eps\cdot \GRAD) \rho_\eps+\rho_\eps
\DIV u_\eps=0,\quad x\in \R^2,t>0\\
\rho_\eps(t=0)=\rho^0,\quad x\in \R^2,
\end{aligned}
\end{equation}
we can then determine the global asymptotics of the Euler system of
pressureless gases (\ref{Euler}).

\begin{Thm}\label{convergence-rho}
{\sl
Consider a  function $b$ satisfying assumptions $(H0)(H1)$.
Let $\rho_0$  be a nonnegative function in
$W^{s-1,\infty}(\R^2)$, and $u_0$ be  a vector-field in
$W^{s,\infty}(\R^2)$ ($s\geq 1$).
  For all $T\leq T^*$ as in Theorem \ref{existence} and all $\eps\leq
\eps_T$, denote by  $(\rho_\eps,u_\eps)$ the solution of
(\ref{Euler}) in $L^\infty([0,T], W^{s-1,\infty}(\R^2)\times
W^{s,\infty}(\R^2))$. Then,
$$\rho_\eps(t,x)- \rho_0(x) \left(1+ t u_0\cdot \nabla \log(x) \cos
\theta_\eps(t,x)-t u_0^\perp\cdot \nabla \log b(x)
\sin \theta_\eps(t,x)
\right)$$
converges strongly to 0 in $L^\infty([0,T]\times \R^2)$, where the phase
$\theta_\eps$ is defined as previously by~(\ref{theta}).
}
\end{Thm}

\bigskip
Let us comment a little on the proof of those theorems, and give the structure
of the paper.

It is quite clear that energy methods will not enable us to have a good
control on the asymptotics of $(\rho_\eps, u_\eps)$, since as soon as we want
a control  on derivatives of $u_\eps$ unbounded terms will appear. So the
most appropriate way to study System (\ref{Euler}) is to rewrite it using the
charcateristics of the flow and to study those characteristics precisely.

Section 2 is therefore devoted to rewriting  System (\ref{Euler}) in
characteristic form, and in the derivation of a few a priori estimates.

In order
to establish the existence of a solution
$(\rho_\eps, u_\eps)$ to system (\ref{Euler}) on a uniform time interval
$[0,T]$, it is enough to see that the solution is well-defined (and smooth) as
long as the flow generates a diffeomorphism $X_\eps(t,.)$
$${dX_\eps\over dt}(t,x)=u_\eps(t,X_\eps(t,x))$$
hence to prove that the characteristics cannot  cross before time $T$. The
precise  estimates on~$DX_\eps$ leading to Theorem 2 are performed in Section
3, they use in a crucial way some results of non-stationary phase type.

The asymptotic behaviour of $u_\eps(t,X_\eps(t,.))$ and
$\rho_\eps(t,X_\eps(t,.))$ is then simply obtained from the explicit
approximation of the characteristics
$X_\eps$, using Taylor expansions for the various fields. In order to
establish the convergence results stated in Theorems \ref{convergence-u} and
\ref{convergence-rho}, the main difficulty is therefore to get a precise
description of the inverse characteristics $X_\eps^{-1}(t,.)$, which is done
in Section 4.


\section{Appropriate formulation of the system}
\setcounter{equation}{0}
\label{sct:caracteristiques}

As pointed out in the introduction, energy estimates do not seem to be
the right angle of attack for our problem. We shall therefore in this
short section present a new formulation of System~(\ref{Euler}), by
means of characteristics (Paragraph~\ref{sct:traj}). In that way
some a  priori estimates
can be deduced immediately (see Paragraph~\ref{sct:apriori}).

To simplify notation, from now on we shall drop the index~$ \e$ in~$
u_\e$ and simply write~$
u$ (and similarly for any other~$ \e$-dependent function).

\subsection{Trajectories associated with the flow }
\label{sct:traj}
Let us write System~(\ref{Euler}) in the
following form:
\begin{equation}
\label{eq:transporteq}
\displaystyle
\begin{array}{c}
\displaystyle\frac{dX}{dt}  =  u(t,X), \quad
X_{|t = 0}  =   x \\
\displaystyle  \frac{d}{dt} (\rho (t,X)) + \rho \nabla \cdot u(t,X) =
 0, \quad \rho_{|t = 0} = \rho_0 \\
 \displaystyle \frac{d}{dt}(u(t,X)) + \frac{b(X)}{\e} u^\perp (t,X)   =
  0, \quad u_{|t = 0} = u_0.
\end{array}
\end{equation}
As seen in Theorem~\ref{thm:existenceeps}, there is a solution to
System~(\ref{Euler})
for a time depending on~$ \e$, and as long as the trajectories do not
intersect we can write in particular
\begin{equation}
\label{eq:formulau}
u(t,X(t,x)) = u_0(x) \cos \left( \frac{\phi (t,x)}{\e}
\right) -  u_0^\perp(x)  \sin \left( \frac{\phi (t,x)}{\e}
\right) ,
\end{equation}
where we have
defined the functions$$ \phi(t,x)=\int_0^t \beta(s,x) \:
ds,\quad \beta (t,x)
\eqdefa b(X(t,x))$$
(these functions are well defined as long as
the characteristics do not cross each other).

If~$ u $ is smooth enough, then~$ \rho $ is uniquely defined by the
transport equation it satisfies. So from now on we can concentrate
on~$ u$ (and~$ X$).
As one of the aims of this article is to prove
Theorem~\ref{existence} (which will be achieved in the next section),
we shall from
now on call~$ T^\e$ the  largest time before which no
characterestic intersect; one of our goals is to prove
that~$ T^\e$ is uniformly bounded from below as $\eps$ goes
to zero.

In the next paragraph we are going to derive
from~(\ref{eq:transporteq}) and~(\ref{eq:formulau}) some easy a priori
estimates for times~$0 \leq t \leq T^\e $, which will help us prove
 Theorem~\ref{existence}
 in the following Section~\ref{sct:X}, and
Theorems~\ref{convergence-u}
 and~\ref{convergence-rho} in Section~\ref{sct:proofconvergence}.

\subsection{A priori estimates}
\label{sct:apriori}
Formula~(\ref{eq:formulau}) immediately enables us to deduce the
following a priori estimate:
\begin{equation}
\label{eq:estiuLinfty}
\|u\|_{L^\infty([0,T^\e[ \times \R^2)} \leq 2 \|u_0\|_{L^\infty} \: ,
\end{equation}
which implies that
\begin{equation}
\label{eq:estidtX}
\displaystyle
\left\|\frac{dX}{dt}\right\|_{L^\infty([0,T^\e[  \times \R^2)} \leq 2
\|u_0\|_{L^\infty}.
\end{equation}
In particular~$ X - x$
remains bounded in space for all times~$0 \leq t < T^\e $,  and we have
\begin{equation}
\label{eq:boundX-xLinfty}
\forall t \in [0,T^\e[, \quad \|X(t,\cdot) - x\|_{L^\infty(\R^2)} \leq 2 t
\|u_0\|_{ L^\infty(\R^2)}.
\end{equation}
Since~$ \beta (t,x) = b(X(t,x))$, we have
\begin{eqnarray}
\displaystyle \|\partial_t \beta \|_{L^\infty([0,T^\e[  \times \R^2)}
&  \leq & \|\nabla
b\|_{L^\infty(\R^2)}  \left\|\frac{dX}{dt}\right\|_{L^\infty([0,T^\e[
  \times  \R^2)}
\nonumber \\
\label{eq:estidtbeta}
\displaystyle
 & \leq  &  2  \|\nabla b\|_{L^\infty(\R^2)}  \|u_0\|_{L^\infty},
\end{eqnarray}
as well as
\begin{equation}
\label{eq:boundsbeta}
\forall t \in [0,T^\e[, \quad \forall x \in \R^2, \quad b_- \leq
\beta(t,x)  \leq \|b\|_{L^\infty(\R^2)}
\end{equation}
with~$b_- $ defined in~$ (H1)$.

Now we are going to look for  an approximation of~$ X$: integrating
 formula~(\ref{eq:formulau}) in time yields
\begin{equation}
\label{eq:defX}
X(t,x) = x +  u_0(x)\int_0^t \cos \left( \frac{\phi
(s,x)}{\e} \right) \: ds-  u_0^\perp(x) \int_0^t  \sin
\left( \frac{\phi (s,x)}{\e}
\right) \: ds
\end{equation}
recalling that~$\displaystyle \phi(t,x) =\int _0^t b(X(s,x))\: ds$.
The following section will be devoted to a  precise study of the
trajectories~$ X$, which will enable us to
infer Theorem~\ref{existence}.


\section{Study of the trajectories}\label{sct:X}
\setcounter{equation}{0}

Formulation~(\ref{eq:transporteq}) of System  (\ref{Burgers})
shows that the study of the Euler system of pressureless
gases with magnetic penalization comes down to a precise
analysis of the characteristics, and in particular of their
invertibility.

In this section we will establish that the trajectories
defined by (\ref{eq:defX}) are invertible on a time interval
$[0,T^\eps[$ with 
$$
\lim_{\eps \to 0} T^\eps = T^*>0,
$$
where~$T^*$ depends on the magnetic field $b$ and on the initial
velocity field $u_0$. This result is based on an asymptotic
expansion of the Jacobian
$$
J(t,x) \eqdefa |\det (DX(t,x)) \,|,
$$
which implies that
$$
\forall t\in [0,T^*[, \quad \liminf_{\eps \to 0} J(t,x) >0.
$$

The asymptotic expansions of $X$ and $DX$ (Paragraphs 3.2
and 3.4)  are obtained  using some
results of non-stationary phase type and the
$L^\infty$-bounds established in Paragraphs 3.1 and~3.3.

\subsection{Bounds on $X(t,\cdot)$}

The first step of the analysis
consists in showing that for any point~$x \in \R^2$, the
characteristic stemming from~$ x$ stays in a ball of size~$ O(\e)$
around~$ x$.  This shows that the rotation has a drastic influence over
the transport by~$ u$.

We have the following proposition.
\begin{Prop}
\label{X-xOe}
{\sl
  Let~$ x \in \R^2$ be given, and let~$ X(\cdot,x)$ be the trajectory
  starting from~$ x$ at time~$ 0$, defined by~(\ref{eq:defX}). As long
  as it  is defined, it
  satisfies
$$
\forall t < \min(T,T^\e), \quad
|X(t,x) - x|  \leq 4 \frac{\e}{b_-}
\|u_0\|_{L^\infty} (1 +  T\frac{\|\nabla
  b\|_{L^\infty}\|u_0\|_{L^\infty}}{b_-})\cdotp
$$
}
\end{Prop}

{\bf Proof of Proposition~\ref{X-xOe}. }  The proof is an immediate
application of the non-stationary phase theorem. As we will be using
such arguments many times in the following, let us state and prove the
following lemma, which will be invoked systematically in the next
sections.
\begin{Lem}
\label{phasenonstat}
{\sl
Let~$ T$ be a given real number, possibly depending on~$ \e$.
Let~$ F $ be a function uniformly bounded in~$ W^{1,\infty}
([0,T],L^\infty(\R^2))$,
and let~$ \beta $  be a positive function, also  uniformly bounded in~$
W^{1,\infty}
([0,T],L^\infty(\R^2))$, and  bounded by below by~$ b_-$. Then for
all~$t \in  [0,T] $  and all~$ x \in \R^2$, the following
bounds hold:
$$
\longformule{
\left|
\int_0^t
F(s,x) \cos \left(
\int_0^s \frac{\beta(s',x)}{\e} \:ds'
\right) \:ds
\right|
}{
\leq \e \left(\frac{\|F(t,\cdot)\|_{L^\infty(\R^2)}}{b_-} + t
\left\|\partial_s
\frac{F(s,\cdot)}{\beta(s,\cdot)}\right\|_{L^\infty([0,t] \times
\R^2)}\right),}
$$
and
$$
\longformule{
\left|
\int_0^t
F(s,x) \sin \left(
\int_0^s \frac{\beta(s',x)}{\e} \:ds'
\right) \:ds
\right|
}{
\leq \e
\left(\frac{\|F(t,\cdot)\|_{L^\infty(\R^2)}+\|F(0,\cdot)\|_{L^\infty(\R^2)}}{b_-
} + t
\left\|\partial_s \frac{F(s,\cdot)}{\beta(s,\cdot)}\right\|_{L^\infty([0,t]
\times \R^2)}\right).}
$$
}
\end{Lem}
{\bf Proof of Lemma~\ref{phasenonstat}. } The proof is a simple
application of the nonstationary phase theorem: an integration by
parts leads to
$$
\longformule{
\int_0^t
F(s,x) \cos \left(
\int_0^s \frac{\beta(s',x)}{\e} \:ds'
\right) \:ds = \e \frac{F(t,x) }{\beta(t,x)} \sin\left(\int_0^t
\frac{\beta(s,x)}{\e} \:ds
\right) }{- \e \int_0^t \partial_s \left(\frac{F(s,x)
}{\beta(s,x)}\right) \sin\left(
\int_0^s
\frac{\beta(s',x)}{\e} \:ds' \right) \: ds,}
$$
and similarly
$$
\longformule{
\int_0^t F(s,x) \sin \left(
\int_0^s \frac{\beta(s',x)}{\e} \:ds'
\right) \:ds = \e  \frac{F(0,x) }{\beta(0,x)} - \e \frac{F(t,x)
}{\beta(t,x)} \cos\left(\int_0^t
\frac{\beta(s,x)}{\e} \:ds
\right)
}{
+ \e \int_0^t \partial_s \left(\frac{F(s,x)
}{\beta(s,x)}\right) \cos\left(
\int_0^s
\frac{\beta(s',x)}{\e} \:ds' \right) \: ds.}
$$
The result follows immediately.

Now let us go back to the proof of Proposition~\ref{X-xOe}. Recalling
formula~(\ref{eq:defX}), we simply apply Lemma~\ref{phasenonstat} to the
case~$ F(t,x) = u_0(x)$ to get
$$
|X(t,x) - x|  \leq 4\e
\frac{\|u_0\|_{L^\infty}}{b_-} + 2 \e t \left\|u_0 \frac{\partial_s
  \beta(s,\cdot)}{\beta^2(s,\cdot)}\right\|_{L^\infty([0,t]\times \R^2)}.
$$
Estimates~(\ref{eq:estidtbeta}) and~(\ref{eq:boundsbeta}) immediately
yield
 Proposition~\ref{X-xOe}.

\subsection{Asymptotics of $X(t,\cdot)$}

The same type of computations based on the non-stationary
phase theorem allows actually to obtain an explicit
approximation of the characteristic $X$ at any order with
respect to $\eps$ (in fact we will stop at order 2 but the
argument can be pushed as far as wanted if necessary).

\begin{Lem}
\label{lem:approxXnextorder}
{\sl
For any point~$ x \in \R^2$ and  any time~$ t \leq
\min(T,T^\eps)$, the following approximation of the
trajectories defined in~(\ref{eq:defX}) holds:
$$
\left| X(t,x) - x - \e \frac{u_0(x)}{b(x)} \sin
\left({\phi(t,x)\over \eps}\right) +
\e
\frac{u_0^\perp(x)}{b(x)}
 \left(1 - \cos\Big({\phi(t,x)\over \eps} \Big)\right)
+ \eps t  v(x) \right|\leq C_T \eps^2,
$$
 where the drift velocity is given by
$$v(x)={1\over 2b^2(x)}\left((u_0^\perp \cdot \nabla b )
u_0(x) - ( u_0
\cdot  \nabla b )  u_0^\perp (x)\right),$$
and~$ C_T$ denotes a
constant depending only on~$ T$, $ u_0$ and~$ b$.
}
\end{Lem}

{\bf Proof of Lemma~\ref{lem:approxXnextorder}.  }
Let us write the following expression for~$ X(t,x)$,
obtained from~(\ref{eq:defX}): we have
$$
X(t,x) = x + R^\e(t,x),
$$
with
\begin{eqnarray}
\label{eq:defRepstx}
R^\e(t,x) & \eqdefa & u_0(x) \int_0^t \cos
\left({\phi(s,x)\over \eps}\right)\:ds - u_0^\perp(x)
\int_0^t \sin \left({\phi(s,x)\over \eps}\right)\:ds\\
 & =  & R^\e_1(t,x) + R^\e_2(t,x) \nonumber.
\end{eqnarray}
We shall only compute the approximation for~$ R^\e_1(t,x)  $, and we
leave~$  R^\e_2(t,x)$ to the reader. By an integration by parts we
have
\begin{equation}
\label{eq:Reps12parts}
R^\e_1(t,x)
=
\e \frac{u_0(x)}{\beta(t,x)}
 \sin \left({\phi(t,x)\over \eps} \right)
- \e u_0(x) \int_0^t \partial_s
 \left(\frac{1}{\beta(s,x)}\right) \sin
\left({\phi(s,x)\over \eps}\right)\:ds.
\end{equation}
The first term is easy to approximate : we have, due to Proposition~\ref{X-xOe}
\begin{equation}
\label{eq:firstpart}
\begin{array}{rl}
\dps \left|\frac{1}{\beta(t,x)} - \frac{1}{b(x)} \right|
&\dps \leq \| \nabla b\|_{L\infty} |X(t,x)-x|\\
&\dps \leq  \frac{4\|\nabla
  b\|_{ L^\infty}\|u_0\|_{ L^\infty}}{b_-^3} (1+T \frac{\|\nabla b\|_{
    L^\infty}  \|u_0\|_{ L^\infty} }{b_-}),
\end{array}
\end{equation}
for all $t\leq \min(T,T^\eps)$ and all $x\in \R^2$.
So~$ \beta (t,x)$ can be replaced by~$ b(x)$ in the first term of~$
R^\e_1$ in~(\ref{eq:Reps12parts}), up to a remainder~$ \e {\mathcal R}^\e$ with~$ \|{\mathcal
  R}^\e\|_{L^\infty([0,T] \times \R^2)} \leq C_T$.

Now we need to approximate the second term. Using the fact that
$$
\partial_s \beta (s,x) = (u \cdot \nabla b)(s,X(s,x))
$$
with
$$
u(s,X(s,x)) = u_0(x) \cos \left({\phi(s,x)\over \eps}\right)
-  u_0^\perp(x)
\sin \left({\phi(s,x)\over \eps}
\right) ,
$$
we can therefore write
$$
- \int_0^t \partial_s
 \left(
\frac{1}{ \beta(s,x)}
\right)
 \sin \left({\phi(s,x)\over \eps}\right)\:ds =  \int_0^t
u_0(x) \cdot \frac{\nabla
 b(X(s,x))}{2b^2(X(s,x))} \sin\left({2\phi(s,x)\over
\eps}\right)
 \: ds
$$
\begin{equation}
\label{eq:decompositionJ}
-  \int_0^t u_0^\perp(x) \cdot  \frac{ \nabla
 b(X(s,x))}{2b^2(X(s,x))} \left(1 -
\cos\left({2\phi(s,x)\over \eps}\right) \right)
 \: ds.
\end{equation}
Note that similar computations  lead to the following
formula, which is useful to estimate $R^\eps_2$:
$$
 \int_0^t \partial_s
 \left(
\frac{1}{ \beta(s,x)}
\right)
 \cos \left({\phi(s,x)\over \eps}\right)\:ds =  \int_0^t
u_0^\perp(x)
 \cdot
\frac{\nabla
 b(X(s,x))}{2b^2(X(s,x))} \sin \left({2\phi(s,x)\over
\eps}\right)
 \: ds
$$
\begin{equation}
\label{eq:decompositionJ2}
-  \int_0^t u_0(x) \cdot  \frac{ \nabla
 b(X(s,x))}{2b^2(X(s,x))} \left(1 +
\cos\left({2\phi(s,x)\over \eps}\right) \right)
 \: ds.
\end{equation}
Both formulas (\ref{eq:decompositionJ}) and
(\ref{eq:decompositionJ2})  show that new harmonics  have
been created by the coupling in the equation.

Let us go back to the estimate of the right-hand side in
(\ref{eq:decompositionJ}).
To estimate the oscillating terms, we use
Lemma~\ref{phasenonstat} with
$$
F_1(s,x) = u_0(x) \cdot \frac{\nabla  b(X(s,x))}{2b^2(X(s,x))} \quad
\mbox{and} \quad  F_2(s,x) = u_0^\perp(x) \cdot \frac{\nabla
  b(X(s,x))}{2b^2( X(s,x))}
\cdotp
$$
We get
$$
\longformule{
\Big| \int_0^t F_1(s,x) \sin\left({2\phi(s,x)\over
\eps}\right)
 \: ds\Big| \leq 2 \e \|u_0\|_{L^\infty} \frac{\|\nabla
   b\|_{L^\infty}}{b_-^3}
}
{
+\e t \|u_0\|_{L^\infty} \left\|\partial_s \frac{\nabla  b
  (X(s,\cdot))}{b^2(X(s,\cdot))}\right\|_{L^\infty([0,t]\times \R^2 )},}
$$
and similarly
$$
\longformule{
\Big| \int_0^t F_2(s,x) \cos\left({2\phi(s,x)\over
\eps}\right)
 \: ds\Big| \leq \e \|u_0\|_{L^\infty} \frac{\|\nabla
   b\|_{L^\infty}}{b_-^3}
}{
+\e t \|u_0\|_{L^\infty} \left\|\partial_s \frac{\nabla  b
  (X(s,\cdot))}{b^2(X(s,\cdot))}\right\|_{L^\infty([0,t]\times \R^2
  )}.
}
$$

By~(\ref{eq:estidtX}) and~(\ref{eq:estidtbeta}) we have
$$
\left\|\partial_s \frac{\nabla  b
  (X(s,\cdot))}{b^2(X(s,\cdot))}\right\|_{L^\infty([0,t]\times \R^2 )}
\leq 2 \|D^2 b\|_{L^\infty} \frac{ \|u_0\|_{L^\infty} }{b_-^2} + 4
\|\nabla b\|^2_{L^\infty}  \frac{ \|u_0\|_{L^\infty} }{b_-^3}\cdotp
$$
Plugging that estimate along with~(\ref{eq:firstpart}) into the
definition  of~$ R^\e_1$ in~(\ref{eq:Reps12parts}), we get finally
$$
R^\e_1 (t,x) =
\e \frac{u_0(x)}{b(x)}
\sin
\left(
\int_0^t
 \frac{\beta(s,x)}{ \e}
\:ds
\right)
 + \e^2 {\mathcal R}^\e(t,x)
- \e u_0(x) \int_0^t \frac{u_0^\perp(x)\cdot \nabla b (X(s,x))}{2 b^2(X(s,x))}
\:ds.
$$
Now we can approximate~$ \displaystyle \frac{\nabla b
(X(s,x))}{b^2
  (X(s,x))}$  by~$ \displaystyle \frac{\nabla b (x)}{
b(x)}$ up
to a  remainder $\eps \RR^\eps$. So we have
$$
- \e u_0(x) \int_0^t \frac{u_0^\perp(x)\cdot \nabla b (X(s,x))}{2 b^2(X(s,x))}
\:ds = -  \e t \frac{u_0^\perp(x)\cdot
  \nabla b (x )}{2 b^2(x)} u_0(x) + \e^2 {\mathcal R}^\e.
$$
The estimate of~$R^\e_2 (t,x) $ is similar and left to the
reader. This ends the proof of Lemma~\ref{lem:approxXnextorder}.

\subsection{A priori estimates on $DX(t,\cdot)$}

A necessary and sufficient condition for $X(t,\cdot)$ to be
invertible is that
$$
J(t,x)\eqdefa |\det (DX(t,x))|
$$
does not cancel. In order to obtain a lower bound on the time
$T^\eps$ (before which the characteristics do no cross each
other), we therefore need to study the behaviour of the 
derivatives~$DX(t,\cdot)$. First of all we derive a uniform $L^\infty$-bound
which will allow to neglect some terms in the asymptotic
expansion.

\begin{Lem}\label{dX-est}
{\sl
  Let~$ x \in \R^2$ be given, and let~$ X(\cdot,x)$ be the trajectory
  starting from~$ x$ at time~$ 0$, defined by~(\ref{eq:defX}). As long
  as it  is defined, it
  satisfies
$$
\forall t < \min(T,T^\e), \quad \forall x\in \R^2,\quad
\|DX(t,x) \|  \leq C_T,
$$
where $C_T$ denotes a constant depending only on $b$, $u_0$ and
$T$.}
\end{Lem}

{\bf Proof of Lemma \ref{dX-est}.}
Differentiating
(\ref{eq:defX}),  leads to
$$
\longformule{ D X(t,x)-Id =   D
u_0(x)\int_0^{t}
\cos
\left({\phi(s,x)\over \eps}
\right) \: ds 
-  D
u_0^\perp(x) \int_0^{t}  \sin
\left({\phi(s,x)\over \eps}
\right) \: ds}{
 -\frac1\eps u_0(x)\int_0^{t} D\phi(s,x)
\sin
\left({\phi(s,x)\over \eps}\right) \: ds  -  \frac1\eps u_0^\perp(x) \int_0^{t}
D\phi(s,x)\cos
\left({\phi(s,x)\over \eps}
\right) \: ds,
}
$$
with
$$
D\phi(s,x)\eqdefa\int_0^s DX(\tau,x)\cdot \nabla b (X(\tau,x))
d\tau.
$$
 Applying the Fubini theorem to both last terms, we can
set this identity in a suitable form to get a Gronwall estimate
\begin{equation}
\label{eq:defdX0}
\begin{array}{rl}
\dps D X(t,x)-Id = &  \dps  D
u_0(x)\int_0^{t}
\cos
\left({\phi(s,x)\over \eps}
\right) \: ds
-  D
u_0^\perp(x) \int_0^{t}  \sin
\left({\phi(s,x)\over \eps}
\right) \: ds\\
& \dps  -\frac1\eps u_0(x)\int_0^{t}
(D X(\tau,x)\cdot
\nabla) b (X(\tau,x))
\int_{\tau}^t \sin
\left({\phi(s,x)\over \eps}
\right) \: ds\: d\tau\\
&\dps -  \frac1\eps u_0^\perp(x) \int_0^{t}
(D X(\tau,x)\cdot
\nabla) b (X(\tau,x))\int_\tau ^t\cos
\left({\phi(s,x)\over \eps}
\right) \: ds \: d\tau.
\end{array}
\end{equation}
From formula (\ref{eq:defX}) we deduce that
$$u_0(x)\int_\tau ^{t} \sin
\left({\phi(s,x)\over \eps}
\right)  \: ds+  u_0^\perp (x) \int_\tau^{t}  \cos
\left({\phi(s,x)\over \eps}
\right)  \: ds= (X(t,x)-X(\tau,x))^\perp.
$$
Plugging this identity back into (\ref{eq:defdX0}) leads to
\begin{equation}
\label{eq:defdX}
\begin{array}{rl}
\dps D X(t,x)-Id = & \dps   D
u_0(x)\int_0^{t}
\cos
\left({\phi(s,x)\over \eps}
\right) \: ds
-  D
u_0^\perp(x) \int_0^{t}  \sin
\left({\phi(s,x)\over \eps}
\right) \: ds\\
&  \dps -\frac1\eps \int_0^{t} (X(t,x)-X(\tau,x))^\perp \otimes
(D X(\tau,x)\cdot
\nabla) b (X(\tau,x))
\: d\tau
\end{array}
\end{equation}
\medskip
As in the  proof of Proposition~\ref{X-xOe}, Lemma~\ref{phasenonstat}
 yields  the
following estimate~: for all $t\leq \min(T,T^\eps)$,
\begin{equation}
\label{reste} \left|D
u_0(x)\int_0^{t}
\cos
\left({\phi(s,x)\over \eps}
\right) \: ds-  D
u_0^\perp(x) \int_0^{t}  \sin
\left({\phi(s,x)\over \eps}
\right) \: ds
\right|\leq C_T\eps
\|\nabla u_0\|_{L^\infty},
\end{equation}
with
$$C_T={4
\over b_-}
\left(1+T{\|\nabla b\|_{L^\infty}
\|u_0\|_{L^\infty} \over b_-}\right).$$

\medskip
From (\ref{eq:defdX}) we then deduce an inequality of
Gronwall type
\begin{equation}
\label{GronwalldX}
\|D X(t,\cdot)\|_{L^\infty}\leq
 1+C_T\eps \|\nabla u_0\|_{L^\infty}+\int_0^{t} \|D
X(\tau,\cdot)\|_{L^\infty}\|\nabla b\|_{L^\infty} \left\|
\frac1\eps (X(\tau,x)-X(t,x))\right\|_{L^\infty} d\tau
\end{equation}
By Proposition \ref{X-xOe},
$$
\forall  t\leq \min(T,T^\eps),\quad \left\|
\frac1\eps (X(t,\cdot)-X(\tau,\cdot))\right\|_{L^\infty} \leq  2C_T
\|u_0\|_{L^\infty},
$$
hence
$$
\|D X(t,\cdot)\|_{L^\infty} \leq (1+C_T\eps \|\nabla
u_0\|_{L^\infty})\exp
\left(2C_T
\|\nabla b\|_{L^\infty}\|u_0\|_{L^\infty} t\right),$$
which is the expected estimate, proving Lemma~\ref{dX-est}.

\subsection{Asymptotics of $DX$ }
In view of the results established in
Lemma \ref{lem:approxXnextorder}, we expect actually
 the derivatives~$\d_i X(t,x)$ to behave asymptotically as
$$\lambda(t,x)+\mu (t,x) \cos \left({\phi(t,x) \over
\eps}\right)+\nu(t,x) \sin \left({\phi(t,x) \over
\eps}\right)
$$
where $\lambda$, $\mu$ and $\nu$ denote some functions which do not depend
on $\eps$.
 Such an
asymptotics can be justified using the same techniques as in the
previous paragraph: let us prove the following lemma.

\begin{Lem}\label{phase-asymp}
{\sl
  Let~$ x \in \R^2$ be given, and let~$ X(\cdot,x)$ be the trajectory
  starting from~$ x$ at time~$ 0$, defined by~(\ref{eq:defX}).
Then, for all $t\leq \min(T^\eps,T)$ and for all $x\in \R^2$,
$$
 {\left\| DX(t,x)-Id -  t u_0\otimes
\nabla
\log b\cos \left({\phi(t,x) \over
\eps}\right)+ t u_0^\perp\otimes
\nabla
\log b  \sin \left({\phi(t,x) \over
\eps}\right)\right\| \leq C_T\eps,}$$
where $C_T$ denotes a constant depending only on $b$, $u_0$ and
$T$.}
\end{Lem}

{\bf Proof of Lemma \ref{phase-asymp}.}
Denote by $g$ the function defined on $[0,T]\times \R^2$  by
$$
g(t,x) \eqdefa DX(t,x)-Id - t u_0\otimes
\nabla
\log b  \cos \left({\phi(t,x) \over
\eps}\right)+t u_0^\perp\otimes
\nabla
\log b  \sin \left({\phi(t,x) \over
\eps}\right)\cdotp 
$$
In view of (\ref{GronwalldX}),  we expect $g$
to satisfy a Gronwall inequality  of the following type
\begin{equation}
\label{Gronwallg}
\|g(t,x)\|\leq  \int_0^{t}\| g(\tau,x)\|
\left\|\frac1\eps(X(t,\cdot)-X(\tau,\cdot))\right\|_{L^\infty}\|\nabla
b\|_{L^\infty}  d\tau +C_T
\eps ,
\end{equation}
for all $t\leq \min(T^\eps,T)$, where $C_T$ denotes a constant
depending only on
$T$,
$u_0$ and
$b$.

\medskip
Let us postpone the proof of this inequality for a while, and
show how it enables us to infer Lemma \ref{phase-asymp}.
 It is easy to see that
$$g(0,x)=0.$$
Applying the Gronwall lemma and using Proposition
\ref{X-xOe} as in (\ref{GronwalldX}) leads  to
$$
\forall t\leq \min(T^\eps,T),\quad \forall x\in \R^2 ,\quad
\|g(t,x)\|\leq C_T
\eps .
$$

\medskip
Now let us go back to the proof of (\ref{GronwalldX}).
We first compute
\begin{equation}
\label{eq:integral}
A(t,x) \eqdefa \int _0^{t} \left(\frac1\eps(X(t,x)-X(\tau,x))\right)^\perp
\otimes \left(g(\tau,x)
\cdot \nabla b(X(\tau,x))\right)d\tau.
\end{equation}
 By Lemma~\ref{lem:approxXnextorder},
$$
\begin{array}{rl}
\dps \frac1\eps(X(t,x)-X(\tau,x))^\perp=&\dps {u_0^\perp\over b}
(x)\left( \sin \left({\phi(t,x) \over
\eps}\right)-\sin \left({\phi(\tau,x) \over
\eps}\right)\right)\\
&\dps -{u_0\over b} (x)\left( \cos \left({\phi(t,x)
\over
\eps}\right)-\cos \left({\phi(\tau,x) \over
\eps}\right)\right)
\\
&\dps -\frac{t-\tau}{2b^2(x)} \left((u_0^\perp\cdot \nabla
 b)u_0^\perp(x)+(u_0\cdot \nabla
 b)u_0(x)\right)+\eps \RR^\eps(t,\tau,x),
\end{array}
$$
where~$  \RR^\eps$ is uniformly bounded in~$ L^\infty([0,T]^2 \times \R^2)$.
Plugging this formula back into the integral~(\ref{eq:integral}) leads to
\begin{equation}
\label{intg}
 \int _0^{t} \left(\frac1\eps(X(t,x)-X(\tau,x))\right)^\perp
\otimes \left(g(\tau,x)
\cdot\nabla b(X(\tau,x))\right)d\tau = A_1(t,x) -  A_2(t,x), 
\end{equation}
with
$$
A_1(t,x) \eqdefa \int
_0^{t}
\left(\frac1\eps(X(t,x)-X(\tau,x))\right)^\perp  \otimes ( D
X(\tau,x) \cdot \nabla b)(X(\tau,x))d\tau 
$$ 
and
$$
\begin{array}{rl}
&\dps  A_2(t,x) \eqdefa  \int_0^{t} \left[ {u_0^\perp\over b}
(x)\left( \sin \left({\phi(t,x) \over
\eps}\right)-\sin \left({\phi(\tau,x) \over
\eps}\right)\right) \right.\\
& \dps-{u_0\over b} (x)\left( \cos \left({\phi(t,x)
\over
\eps}\right)-\cos \left({\phi(\tau,x) \over
\eps}\right)\right) \\
& \dps \left. \qquad -\frac1{2b^2(x)} (t-\tau)\left((u_0^\perp\cdot \nabla
 b)u_0^\perp(x)+(u_0\cdot \nabla
 b)u_0(x)\right)+\eps \RR^\eps(t,\tau,x) \right]\\
&\dps \quad \otimes \left[\nabla b(X(\tau,x)) +(u_0(x)\cdot
\nabla b(X(\tau,x)))
\nabla
\log b (x)\tau \cos \left({\phi(\tau,x) \over
\eps}\right)\right.\\
&\dps \qquad \left. -(u_0^\perp(x)\cdot \nabla b(X(\tau,x)))
\nabla
\log b (x)\tau \sin \left({\phi(\tau,x) \over
\eps}\right) \right]\:d\tau.
\end{array}
$$

From (\ref{eq:defdX}) and (\ref{reste}) we deduce that
$$
A_1(t,x) = Id-DX(t,x)
$$
up to terms of order $\eps$.

In order to estimate the second term~$A_2(t,x)$, we use again a non-stationnary phase
theorem. Since the trajectories lie in balls of size~$\eps$,
$$|\nabla b(X(\tau,x))-\nabla b (x) |\leq C_T
\|D^2b\|_{L^\infty} \eps ,$$ for all $x\in \R^2$ and all $\tau
\leq t\leq \min (T^\eps,T)$. Then, as
$\d_s
\beta$ is uniformly bounded according to~(\ref{eq:estidtbeta}), 
Lemma~\ref{phasenonstat} shows that
$$
\begin{array}{rl}
&A_2(t,x)=\dps \left[t{u_0^\perp\over b}
(x) \sin \left({\phi(t,x) \over
\eps}\right) -t{u_0\over b} (x) \cos \left({\phi(t,x)
\over
\eps}\right)-\frac{t^2}2 v^\perp(x) \right]\otimes \nabla b(x)\\
&\dps +\int_0^t \left({u_0^\perp\over b}(x)\sin
\left({\phi(\tau,x)
\over
\eps}\right)\right) \otimes
\left((\tau u_0^\perp(x)\cdot \nabla b(x))
\nabla
\log b (x) \sin \left({\phi(\tau,x) \over
\eps}\right)\right)d\tau\\
&\dps +\int_0^t \left({u_0\over b}(x)\cos
\left({\phi(\tau,x)
\over
\eps}\right)\right) \otimes
\left((\tau u_0(x)\cdot \nabla b(x))
\nabla
\log b (x) \cos \left({\phi(\tau,x) \over
\eps}\right)\right)d\tau +\eps \RR^\eps (t,x)
\end{array}
$$
Using the identities
$$\cos^2 \phi =\frac12 (1+\cos(2\phi)),\quad \sin^2\phi=\frac12 (1-\cos
(2\phi)),$$
we then obtain that
$$ 
A_2(t,x) = \left( t{u_0^\perp\over b}
(x) \sin \left({\phi(t,x) \over
\eps}\right) -t{u_0\over b} (x) \cos \left({\phi(t,x)
\over
\eps}\right)\right)\otimes \nabla b (x) 
$$
up to terms of order $\eps$.

Then  (\ref{intg}) can be rewritten
$$\int _0^{t}
\left(\frac1\eps(X(t,x)-X(\tau,x))\right)^\perp\otimes
\left( g(\tau,x)\cdot \nabla  b(X(\tau,x))\right) d\tau=-g(t,x)
+\eps
\RR^\eps(t,x),$$ which implies immediately (\ref{Gronwallg}) and yields
Lemma \ref{phase-asymp} as explained above.

\subsection{Existence on a uniform time interval}
As an immediate corollary of Lemma \ref{phase-asymp} we obtain
that $X(t,\cdot)$ is a diffeomophism of $\R^2$ on a uniform time
interval. Indeed, $ DX$  is invertible as long as
$$
\|DX-Id\|_{L^\infty} <1.
$$
\begin{Cor} \label{inversion}
{\sl
Consider a  function $b$ satisfying assumptions $(H0)(H1)$.
Let $(\rho_0,u_0)$ be respectively a nonnegative function
of~$W^{s-1,\infty}(\R^2) $ and a vector-field of~$W^{s,\infty}(\R^2)$
($s\geq 1$).  Then, for all $T<
\|u_0\|_{L^\infty}^{-1} \| \nabla b\|_{L^\infty}^{-1}$, there
exists $\eps_T>0$ such that System~(\ref{Euler}) admits a unique
solution
$(\rho_\eps,u_\e) \in  L^\infty([0,T], W^{s-1,\infty}(\R^2) \times W^{s,\infty}(\R^2))$ for all $\eps \leq \eps_T$.
}
\end{Cor}

{\bf Proof of Corollary \ref{inversion}.}
By Lemma \ref{phase-asymp}, the trajectories  defined
by~(\ref{eq:defX}) are
continuously differentiable and  satisfy for all
$t\leq T$ and  all $x\in
\R^2$,
$$
 {\left\| DX(t,x)-Id - t u_0\otimes
\nabla
\log b \cos \left({\phi(t,x) \over
\eps}\right)+ t u_0^\perp\otimes
\nabla
\log b\sin \left({\phi(t,x) \over
\eps}\right)\right\| \leq C_T\eps.}$$
This implies in particular the following estimate on  the
Jacobian
$J(t,x)\eqdefa |\det(DX(t,x))|$~:
$$
\left| J(t,x)-1- t u_0\cdot
\nabla
\log b\cos \left({\phi(t,x) \over
\eps}\right)+ t u_0^\perp\cdot
\nabla
\log b  \sin \left({\phi(t,x) \over
\eps}\right)\right| \leq C_T\eps.
$$
Then for $T< \|u_0\|_{L^\infty}^{-1} \| \nabla
b\|_{L^\infty}^{-1}$, there exists $\eps_T$ such that
$$\forall \eps \leq \eps_T,\quad \forall t\in [0,T],\quad
\sup_{x\in \R^2} |J(t,x)-1|<1,$$
which means that $X$ is a $C^1$-diffeomorphism of $\R^2$.

\medskip
Moreover, from formula (\ref{eq:defX})
we can deduce by induction that $X(t,.)$ (and consequently its
inverse $X^{-1}(t,\cdot)$) is smooth, its regularity being the same
as the regularity of the initial velocity field $u_0$. Then the
vector field $u$ given by
$$
\longformule{
u(t,x)=u_0(X^{-1}(t,x))\int_0^t \cos
\left({\phi(s,X^{-1}(t,x))\over
\eps }\right)ds}{-u_0^\perp (X^{-1}(t,x))\int_0^t \cos
\left({\phi(s,X^{-1}(t,x))\over
\eps} \right)ds 
}
$$ belongs to $L^\infty([0,T],
W^{s,\infty}(\R^2))$ and it is easy to check that it satisfies
System (\ref{Burgers}) in  strong sense.

\medskip
The density $\rho$ is then obtained as the strong solution of the
linear transport equation
$$\d _t \rho +u\cdot \nabla \rho +\rho \nabla \cdot u =0$$
whose
   coefficients belong to
$L^\infty([0,T],W^{s-1,\infty}(\R^2))$, with initial data in
$W^{s-1,\infty}(\R^2)$. It therefore  stays
in
$L^\infty([0,T],W^{s-1,\infty}(\R^2))$. We emphasize once again
that no uniform bound on $(\rho,u)$ is available in
$L^\infty([0,T],W^{s-1,\infty}(\R^2)\times W^{s,\infty}(\R^2))$.

Theorem~\ref{existence} is proved.
\begin{Rem}
{\sl
The supremum of the life span of the solutions corresponds to a
crossing phenomenon, to be compared with the caustic in
geometrical optics. Beyond this time, the
differential system
$$\left\{
\begin{array}l
\dot X = \xi\\
\dot \xi =\xi \wedge b  \\
\end{array}\right.$$
with initial data $(x,u_0(x))_{x\in \R^2}$ still admits  a
unique smooth  solution, but the application
$(X(t,x),\xi(t,x))\mapsto X(t,x)$ is no longer injective, it
cannot be lifted. The hyperbolic system~(\ref{Burgers})  no
longer has a  solution.
}
\end{Rem}


\section{Study of the asymptotics of~$ u_\e$ and $\rho_\eps$}
\setcounter{equation}{0}
\label{sct:proofconvergence}

Let $T<T^*=
\|u_0\|_{L^\infty}^{-1} \| \nabla b\|_{L^\infty}^{-1}$ be
fixed. Then, for any $\eps \leq \eps_T$ as in Corollary
\ref{inversion}, the solution $(\rho,u)$ of System
(\ref{Euler}) with initial data $(\rho_0,u_0)\in
W^{s-1,\infty}(\R^2) \times W^{s,\infty}(\R^2)$ belongs to
$L^\infty([0,T], W^{s-1,\infty}(\R^2) \times
W^{s,\infty}(\R^2))$. Then it makes sense to study their
asymptotic behaviour as $\eps \to 0$, and the aim of this section is
to  prove
Theorems~\ref{convergence-u} 
and~\ref{convergence-rho}.

Paragraph~\ref{sct:asymptotuX} is devoted to the asymptotics
of~$u(t,X(t,x)) $ ad~$\rho (t,X(t,x))$. The last paragraph 
consists in inverting the characteristics in order to infer
Theorems~\ref{convergence-u}  
and~\ref{convergence-rho}.

\subsection{Asymptotics of~$ u(t,X)$ and $\rho(t,X)$}
\label{sct:asymptotuX}

From the characteristic formulation of System~(\ref{Euler})
and the asymptotic expansion of $X(t,\cdot)$ we immediately
deduce the asymptotic behaviour of $u(t,X(t,\cdot))$ and
$\rho(t,X(t,\cdot))$.

\begin{Prop}
\label{pro:convergenceuX}
{\sl
Consider a  function $b$ satisfying assumptions $(H0)(H1)$.
Let $u_0$ be  a vector-field in $W^{s,\infty}(\R^2)$
($s\geq 1$). For all~$ T<T^*$ and~$ \e \leq \e_T$ as in
Theorem~\ref{existence}, denote by  $u$ the solution
of~(\ref{Burgers})
in~$L^\infty([0,T],
W^{s,\infty}(\R^2))$. Then
$$
u(t,X(t,x))- \left( u_0(x) \cos (\widetilde \phi_\eps(t,x))-u_0^\perp(x) \sin
(\widetilde \phi_\eps(t,x))
\right)
$$
converges strongly to 0 in $L^\infty([0,T]\times \R^2)$, at speed~$ O(\e)$, where the phase
$\widetilde \phi_\eps$ is defined by
\begin{equation}
\label{phitilde}
\widetilde \phi_\eps(t,x)={b(x)t\over \eps}-t
\left(u_0^\perp(x)\cdot
\GRAD\right)   \log b(x)
.
\end{equation}
}
\end{Prop}

{\bf Proof of Proposition~\ref{pro:convergenceuX}. }
Let us first recall that
$$u(t,X(t,x))=u_0(x)\int_0^t \cos
\left({\phi(s,x)\over
\eps }\right)ds-u_0^\perp (x)\int_0^t \sin
\left({\phi(s,x)\over
\eps} \right)ds ,$$
where the phase $\phi$ is given by
$$\phi(t,x)=\int_0^t b(X(s,x)) ds .$$
Then in order to establish Proposition
\ref{pro:convergenceuX}, we have to
  approximate the phase. By
 Lemma~\ref{lem:approxXnextorder},
$$
\longformule{
b(X(t,x)) = b(x) + \frac{\e u_0(x)}{b(x)}
  \sin \left({\phi(t,x)\over \eps}\right) \cdot \nabla b(x)}{
- \frac{\e u_0^\perp(x)}{b(x)} \left(1 -  \cos
\left({\phi(t,x)\over \eps}\right)\right) \cdot \nabla b(x)
+ \e^2   {\mathcal R}^\e(t,x), }
$$
noticing that~$ v \cdot\nabla b = 0$.
It follows that
$$
\longformule{
u(t,X(t,x)) =
u_0(x) \cos
\Big(\frac{b(x)t}{\e} +  \frac{
  u_0(x) \cdot \nabla b(x)}{b(x)}
\int_0^t  \sin \left({\phi(s,x)\over \eps}
\right) \:ds
}{
-  \frac{
  u_0^\perp(x) \cdot \nabla b(x)}{b(x)} \int_0^t \left(1 -
 \cos \left({\phi(s,x)\over \eps}\right) \right)\:ds
 + \e  {\mathcal R}^\e(t,x)\Big)
}
$$
$$
\longformule{
-  u_0^\perp(x) \sin\Big(\frac{b(x)t}{\e} +  \frac{
  u_0(x) \cdot \nabla b(x)}{b(x)} \int_0^t  \sin
\left({\phi(s,x)\over \eps}\right) \:ds }{
-  \frac{
  u_0^\perp(x) \cdot \nabla b(x)}{b(x)} \int_0^t \left(1
-  \cos \left({\phi(s,x)\over \eps}\right) \right)\:ds
+ \e  {\mathcal R}^\e(t,x)\Big)
}
$$
Finally remembering that due to Lemma~\ref{phasenonstat}
$$
\left|\int_0^t  \sin \left({\phi(s,x)\over \eps}\right) \: ds
\right| + \left|\int_0^t  \cos \left({\phi(s,x)\over
\eps}\right)  \: ds \right| \leq \e {\mathcal R}^\e (t,x),
$$
with the usual uniform bounds on~$ {\mathcal R}^\e$, yields
Proposition~\ref{pro:convergenceuX}.

\bigskip
The asymptotic behaviour of $\rho$ is obtained in a similar
way using the fact that $\rho$ is proportionnal to the
Jacobian $J(t,x)=|\det DX(t,x)|$.

\begin{Prop}
\label{pro:convergencerX}
{\sl
Consider a  function $b$ satisfying assumptions $(H0)(H1)$.
Let $\rho_0$  be a nonnegative function in
$W^{s-1,\infty}(\R^2)$, and $u_0$ be  a vector-field in
$W^{s,\infty}(\R^2)$ ($s\geq 1$). For all~$ T<T^*$ and~$ \e
\leq \e_T$ as in Theorem~\ref{existence}, denote by
$(\rho,u)$ the solution of~(\ref{Euler})
in~$L^\infty([0,T],
W^{{s-1},\infty}(\R^2))~$ and~$L^\infty([0,T],
W^{s,\infty}(\R^2))$ respectively. Then
$$
\rho(t,X(t,x))- \rho_0(x) \left( 1+ tu_0\cdot \nabla \log
b(x)
\sin (\widetilde
\phi_\eps(t,x))-t u_0^\perp\cdot \nabla \log b(x) \cos (\widetilde
\phi_\eps(t,x))
\right)
$$
converges strongly to 0 in $L^\infty([0,T]\times \R^2)$, where the phase
$\widetilde \phi_\eps$ is defined as previously by~(\ref{phitilde}).}
\end{Prop}

{\bf Proof of Proposition~\ref{pro:convergencerX}. }
As long as the solution of (\ref{Euler}) is regular, the
equation governing~$\rho$ can be rewritten
$${d\over dt} (\log \rho)=\nabla\cdot u,$$
where ${d\over dt}$ denotes as usual the derivative along
the trajectories associated with the flow. Of course, the
Liouville theorem implies that the equation on the Jacobian
of the flow states
$${d\over dt} (J)=\nabla\cdot u.$$
Then, for all $\eps \leq \eps_T$, all $t\in [0,T]$ and all
$x\in \R^2$,
$$\rho(t,X(t,x))=\rho_0(x) J(t,x),$$
since $J_0(t,x)=\det (Id)=1$.

From Lemma \ref{phase-asymp} we then deduce that
\begin{equation}
\label{rho-asymp}
\left|\rho(t,X(t,x))-\rho_0(x)\left(  1+
tu_0\cdot
\nabla
\log b(x)
\sin \left({\phi(t,x)\over \eps}\right)-t u_0^\perp\cdot
\nabla \log b(x)
\left({\phi(t,x)\over \eps}\right)
\right) \right|\leq C_T\eps.
\end{equation}
Plugging the approximation of the phase obtained previously
\begin{equation}
\label{eq:moreprecise}
{\phi(t,x)\over \eps}=\widetilde \phi_\eps+\eps
\RR^\eps(t,x)
\end{equation}
back into  formula (\ref{rho-asymp}) leads then to the
expected asymptotics.

\subsection{Inversion of the characteristics}
In this section we shall prove Theorems~\ref{convergence-u}
and~\ref{convergence-rho}.   From now on~$ T^*$ is the time given by
Theorem~\ref{existence}, and we will call~$ T$ any time smaller than~$
T^*$ (in the following we will also suppose~$ \e \leq \e_T$ as given
in Theorem~\ref{existence}).

Let~$ X^{-1}(t,x)$ be the point at time~$ 0$ of the trajectory
reaching~$ x$ at time~$ t$.  By Proposition~\ref{pro:convergenceuX}, we have
$$
u(t,x) = u_0 ( X^{-1}(t,x))  \cos \left(\widetilde
\phi_\eps(t, X^{-1}(t,x))\right)-u_0^\perp(  X^{-1}(t,x)) \sin
\left(\widetilde \phi_\eps(t, X^{-1}(t,x))\right) + \e {\mathcal R}^\e(t,x)
$$
with the usual uniform bounds on~$  {\mathcal R}^\e$. That remainder 
function~$ {\mathcal R}^\e(t,x)$ is liable to change from line to line
in this paragraph.

By Proposition~\ref{X-xOe}  there is a constant~$ C_T$ (depending
on~$ T$, $ u_0$ and~$ b$), such that
\begin{equation}
\label{eq:X-1closetox}
\forall x \in \R^2, \quad \forall t \in [0,T] , \quad | X^{-1}(t,x) -
x| \leq C_T
\e,
\end{equation}
so
we can write
rather
$$
u(t,x) = u_0 (x)  \cos \left(\widetilde
\phi_\eps(t, X^{-1}(t,x))\right)-u_0^\perp(x) \sin
\left(\widetilde \phi_\eps(t, X^{-1}(t,x))\right) + \e {\mathcal R}^\e(t,x).
$$
By definition of~$ \widetilde \phi_\e$ in~(\ref{phitilde}),
 we have, using again~(\ref{eq:X-1closetox}),
$$
\widetilde \phi_\eps(t, X^{-1}(t,x)) = \frac{b(x)t}{\e} + \frac{t}{\e}
(X^{-1}(t,x) - x) \cdot \nabla b(x) - t u_0^\perp (x) \cdot \nabla
\log b(x) +  \e {\mathcal R}^\e(t,x)
$$
hence defining
$$
 \widetilde \theta_\e (t,x) \eqdefa  \frac{t}{\e}
(X^{-1}(t,x) - x) \cdot \nabla b(x) - t u_0^\perp (x) \cdot \nabla
\log b(x) 
$$
we have
\begin{equation}
\label{eq:phitheta}
\widetilde \phi_\eps(t, X^{-1}(t,x)) = \frac{b(x)t}{\e} + \widetilde \theta_\e
(t,x)   +  \e {\mathcal R}^\e(t,x).
\end{equation}
Now we shall try to make~$ \widetilde  \theta_\e$ more precise.
According to Lemma~\ref{lem:approxXnextorder} and the approximation
for the phase derived in the previous paragraph, we have
$$
\longformule{
x -  X^{-1}(t,x) = \e \frac{u_0(x)}{b(x)} \sin \left(\widetilde
\phi_\eps(t, X^{-1}(t,x))  + \e{\mathcal R}^\e(t,x) \right)
}{ - \e
\frac{u_0^\perp(x)}{b(x)} \left(1 - \cos \left(\widetilde
\phi_\eps(t, X^{-1}(t,x))  + \e{\mathcal R}^\e(t,x) \right)\right)- \e t v(x) + \e^2 {\mathcal R}^\e(t,x)
}
$$
where again we have used~(\ref{eq:X-1closetox}).
So we obtain, using the fact that~$v \cdot \nabla b = 0$, 
$$
 \widetilde  \theta_\e (t,x) = -t u_0(x) \cdot \nabla 
\log b(x)  \sin \left(\widetilde
\phi_\eps(t, X^{-1}(t,x))  + \e{\mathcal R}^\e(t,x) \right) 
$$
$$
+ t u_0^\perp(x) \cdot \nabla 
\log b(x)  \cos \left(\widetilde
\phi_\eps(t, X^{-1}(t,x))  + \e{\mathcal R}^\e(t,x) \right) + \e {\mathcal R}^\e(t,x)
$$
which by~(\ref{eq:phitheta}) yields directly the
result~(\ref{theta}), defining~$
 \displaystyle \theta_\e (t,x) 
= \widetilde  \theta_\e (t,x) +\frac{b(x)t }{\e}\cdotp
$ 

Theorem~\ref{convergence-u} is proved.

The proof of 
Theorem~\ref{convergence-rho} is now immediate: we use the formula
obtained in
Proposition~\ref{pro:convergencerX} and replace~$\rho_\e(t,X(t,x))$
by~$\rho_\e(t,x)$ using the above formulation of~$
X^{-1}(t,x)$. The result follows.

\end{document}